\documentclass[11pt]{amsart}
\usepackage{graphicx,amssymb,mathrsfs,amsmath}
\usepackage{stmaryrd}

\usepackage{mathrsfs}
\usepackage{amscd}
\usepackage{amscd, amssymb, amsmath, amsthm, graphics}
\usepackage{amsmath,amsfonts,amsthm,amssymb}
\usepackage{latexsym,amsmath}
\usepackage{graphicx,psfrag}
\usepackage[all]{xy}
\usepackage{latexsym}
\usepackage{enumerate}

\newtheorem{theorem}{Theorem}[section]
\newtheorem{lemma}[theorem]{Lemma}
\newtheorem{prop}[theorem]{Proposition}
\newtheorem{cor}[theorem]{Corollary}

\theoremstyle{definition}
\newtheorem{definition}[theorem]{Definition}
\newtheorem*{example}{Example}

\newtheorem{remark}[theorem]{Remark}

\newtheorem*{claim}{Claim}

\numberwithin{equation}{section}

\newcommand{\CC}{\mathbb{C}}
\newcommand{\RR}{\mathbb{R}}

\newcommand{\ZZ}{\mathbb{Z}}

\newcommand{\MM}{{\bf M}}
\newcommand{\JJ}{\mathscr{J}}
\newcommand{\Jst}{\mathscr{J}_{\beta,c}}
\newcommand{\thmone}{\begin{theorem}}
\newcommand{\thmtwo}{\end{theorem}}
\newcommand{\lemmaone}{\begin{lemma}}
\newcommand{\lemmatwo}{\end{lemma}}
\newcommand{\defone}{\begin{definition}}
\newcommand{\deftwo}{\end{definition}}
\newcommand{\corone}{\begin{cor}}
\newcommand{\cortwo}{\end{cor}}
\newcommand{\pfone}{\begin{proof}}
\newcommand{\pftwo}{\end{proof}}
\newcommand{\rkone}{\begin{remark}}
\newcommand{\rktwo}{\end{remark}}

\newcommand{\onehalf}{\left(\begin{array}{cc}}
\newcommand{\theother}{\end{array}\right)}
\newcommand{\onedis}{\begin{displaymath}}
\newcommand{\twodis}{\end{displaymath}}
\newcommand{\oneeq}{\begin{equation}}
\newcommand{\twoeq}{\end{equation}}
\newcommand{\nono}{\noindent}

\newcommand{\bpl}{\bar{\partial}}
\newcommand{\pl}{\partial}
\newcommand{\sq}{\subseteq}

\begin{document}
\abovedisplayskip=6pt plus 1pt minus 1pt \belowdisplayskip=6pt plus
1pt minus 1pt
\thispagestyle{empty} \vspace*{-1.0truecm} \noindent
\vskip 9mm

\begin{center}{\large\bf Note on a theorem of Bangert

\footnotetext{\footnotesize $^*$ The first author is supported by NSF grant.}} \end{center}

\vskip 5mm
\begin{center}{\bf $^*$Tian-Jun Li\\
{\small\it Department of Mathematics, University of Minnesota, Twin
Cities,
MN {\rm 55455}\\
E-mail: tjli@umn.edu\quad }}\end{center}

\begin{center}{\bf Weiwei Wu\\
{\small\it Department of Mathematics, University of Minnesota, Twin
Cities,
MN {\rm 55455}\\
E-mail: weiwei@umn.edu\quad }}\end{center}

\vskip 1 mm

\noindent{\small {\small\bf Abstract:} We generalize Bangert's
non-hyperbolicity result for uniformly tamed almost complex
structures on standard symplectic $\RR^{2n}$ to asymtotically standard
symplectic manifolds.  \ \

\vspace{1mm}\baselineskip 10pt


\noindent{\small\bf MR(2000) Subject Classification:}  Primary:
58D10 Secondary: 32Q45, 53D45 \ \ {\rm }}

\section{Introduction}

Let $M$ be a smooth manifold and $J$ an almost complex structure on $M$.  Any given $v\in TM$ is tangent
to a $J$-holomorphic disc (see \cite{IR}, \cite{NW}).  Therefore we could define
Kobayashi-Royden pseudo-norm $||\cdot||_{k,J}$ of tangent vectors as:

$$||v||_{k,J}=inf\{R^{-1}:\exists \text{a pseudo-holomorphic } f:B_R\rightarrow(M,J), df(\partial_x)=v\}$$

\nono Here $B_R\subset \CC$ is the ball of radius $R$, $\partial_x$ the
unit tangent vector along the real axis.  This in turn induces a
pseudo-distance on $M$, which is exactly the Kobayashi
(pseudo-)distance.  If such an induced distance is a metric, then $J$ is
called (Kobayashi-)hyperbolic.  Aother well-known notion of complex hyperbolicity
is \textit{Brody hyperbolicity}, which means the absence of holomorphic
lines in $(M,J)$.  These two notions are equivalent on compact manifolds, see \cite{Lang}.

In a non-compact Riemannian manifold $(M,g)$, it seems natural to
consider the uniform property of hyperbolicity.  We say $J$ is
\textit{($g$-)uniformly hyperbolic}, if $||v||_g\leq C\cdot||v||_{k,J}$,
for all $v\in TM$ and a constant $C$ independent of $v$.  This is equivalent to saying that the set

$$\{R\text{ }|\exists \text{a pseudo-holomorphic map } f:B_R\rightarrow(M,J),\text{ }||f'(0)||_g=1\}$$

\nono is bounded.  In particular, it is easy to see that for a
non-compact complete manifold, the completeness of Kobayashi distance (see \cite{IR}, \cite{DI}, \cite{GS})
 follows from the uniform hyperbolicity.  If $M$ is non-compact and $(M,J)$ is not Brody hyperbolic, that is,
there exist a complex line, then for any Reimannian metric $g$, $J$ cannot be uniformly hyperbolic.

The notion of uniform hyperbolicity was studied first by V.Bangert for
the class of $J$ uniformly tamed by standard symplectic form $\omega_0$ on $(\RR^{2n},g_{eud})$ in \cite{Bangert},
where he proved that none of them is uniformly hyperbolic.  Recall
that when $(M, \omega)$ is a symplectic manifold, an almost-complex
structure $J$ on $M$ satisfying $\omega(v,Jv)>0$ for any tangent vector
$v\in TM$ is called $\omega$-\textit{tamed} (see \cite{Gromov}),
while the notion \textit{uniformly tamed} introduced by Bangert in \cite{Bangert} is
more suitable in the setting that $(M,\omega)$ is a non-compact symplectic manifold endowed
with a preferred metric $g$ (see Definition \ref{Uniformlytamedpair}).
Motivated by \cite{Bangert}, in this note we generalize
Bangert's result to asymptotically standard symplectic manifolds.

\defone\label{noncompact} A (non-compact) symplectic manifold with Riemannian metric $(M, g, \omega)$
is \textit{asymptotically standard} if for some compact subset $K$ in
$M$, there is a map $\phi_M:(M\backslash K,g,\omega)\rightarrow(\RR^{2n}\backslash B_{\bar{R}},g_{eud},\omega_0)$ which
 is both a symplectomorphism and an isometry.\deftwo

 It should be noted that in the case when $(M,\omega)$ is symplectic aspherical, it is proved by
Eliashberg, Floer and McDuff that $M$ is diffeomorphic to $\RR^{2n}$
(cf. \cite{M},\cite{MS2}).
 Our main result in this paper reads as follows:

\thmone\label{main} For an asymptotically standard symplectic manifold $(M,
g, \omega)$, any uniformly tamed almost complex structure is not
uniformly hyperbolic. \thmtwo


In essence what we need is to construct a sequence of holomorphic disks with radius approaching infinity,
while the norm of the differential at the origin remains $1$.
The starting point of the proof is  the ``almost-K\"ahler cut"
technique due to D. Burns, V. Guillemin, and E. Lerman in
\cite{BGL}, which provides a convenient way
 to compactify the manifold while keeping track of the almost complex structure.
By  \cite{HR} (see also \cite{Mcduff2}), the compactification is
easily seen to be rationally connected.  We thus obtain a sequence
of holomorphic disks inside larger and larger parts of $M$.
 One then uses the reparametrization process developed in \cite{Bangert} to obtain needed estimates for derivatives.

There are a few remarks we need to make about Theorem \ref{main}.
  Similar problems have been explored by
various authors, see \cite{Bio}, \cite{He}. Especially, when
$(\omega, J)$ is a compatible pair, uniform hyperbolicity is
equivalent to the notion of \textit{almost-K\"ahler hyperbolicity}
in \cite{Bio} with $g$ being the induced metric by the $(\omega,
J)$. A notion of Floer theoretical symplectic hyperbolicity was
introduced in \cite{Bio}, which is quite different, though related.
 Bangert's result is then extended by
Biolley to many generalized Stein manifolds in the class of
compatible almost-complex structures, with some restrictions on
capacity.


It is notable that similar patterns appear in the far-reaching
symplectic field theory, see \cite{symfield},
\cite{compactness}, \cite{foliation}, \cite{Hind}, \cite{Wel}, etc.
The ``almost K\"ahler cut" technique we employ in this paper should
be viewed somehow as the process of ``stretching the neck" in
symplectic field theory, with the addition feature of retaining the almost complex
structure.

\textit{Acknowledgement:} The authors would like to thank
Shinichiroh Matsuo for suggesting the paper of V. Bangert, and
helpful discussions about results in his upcoming papers on
analogous problems in gauge theory settings.  Gratitude is also due
to Jianxun Hu and Guangcun Lu for inspiring conversations.  Special thanks
to the anonymous referees who gave numerous suggestions making the exposition
of this note much more clear. The second author would also like to
thank Josef Dorfmeister, Guoyi Xu and Weiyi Zhang for their interests in this work.

\section{Uniformly tamed almost complex structures}
\defone\label{Uniformlytamedpair} \textup{ Let $(M,g)$ be a Riemannian manifold equipped with
a symplectic form $\omega$, and an almost complex structure $J$ over
$M$.  We say $J$ and $\omega$ are \textit{($g$-)uniformly tamed} by each
other if}:

\begin{enumerate}[(1)]

\item (uniform boundedness of $\omega$) For some $\alpha>0$, $\omega(u,v)\leq\alpha||u||_g||v||_g$

\item (uniform boundedness of $J$) For some $\beta>0$,   $||Jv||_g\leq \beta||v||_g$

\item (uniform tameness) For some $c>0$, $\omega(v,Jv)\geq c||v||_g^2$
\end{enumerate}

\deftwo

\begin{remark}\label{Theremark}
From the definition above, one easily deduce that in the case of
uniformly tamed pair, the metric defined by
$\tilde{g}(u,v)=\frac{1}{2}(\omega(u,Jv)-\omega(Ju,v))$ is
equivalent to the metric $g$, in that there is a constant $a>0$
satisfying $\frac{1}{a}||v||_{\tilde{g}}\leq||v||_g\leq
a||v||_{\tilde{g}}$ for all tangent vector $v$.  Moreover, a
$J$-holomorphic curve has the same symplectic area as its associated
$\tilde{g}$-area.  This follows from that the pull-back metric has
volume form $\omega((f_R)_*(v),J(f_R)_*(v))ds\wedge dt$, where
$ds,dt$ are dual to $v,iv$, respectively (see also the proof of
Lemma 1.2 in \cite{Bangert}).  In what follows we will call $\beta$
the \textit{bounding constant} and $c$ the \textit{taming constant}.
\end{remark}

Now we focus on almost-complex structures on $\RR^{2n}$.  Let
$\JJ(\omega_0)$ denote the space of complex structures tamed by
$\omega_0$, and $\mathscr{J}_{\beta,c}$
those uniformly tamed by $\omega_0$ with bounding constant $\beta$
and taming constant $c$. Following \cite{Bangert}, note that $\Jst$
is closed and bounded as a subset of $\mathscr{J}(\omega_0)$,
therefore compact. Since the tangent bundle of $\RR^{2n}$ is
trivial, a uniformly tamed almost-complex structure $J_{\beta, c}$
with taming constant $\beta$ and bounding constant $c$ on a subset
$S$ of $\RR^{2n}$ yields a map from $S$ to
$\mathscr{J}_{\beta,c}$. Conversely, a compact subset of
$\mathscr{J}(\omega_0)$ is bounded with a uniformly tamed constant.
Therefore,

\lemmaone\label{uniformlytamenessandcompactness} A uniformly
tamed almost-complex structure $J_{\beta, c}$ over
$S\subset\RR^{2n}$ is equivalent to a map with precompact image
$f:S\rightarrow \mathscr{J}_{\beta,c}$. \lemmatwo

\section{A variation of almost K\"ahler  cut}

\subsection{Almost K\"ahler cut}
We now give a brief review of the technique of almost-K\"ahler cuts.
Readers are referred to \cite{BGL} for a more detailed
exposition on the subject.  It should be noted that although
\cite{BGL} concerns with K\"ahler manifolds only, the results we use
in this note can be easily (even directly) adapted
 to the context of almost-K\"ahler  manifolds.
\\

We first review the almost-K\"ahler  reduction process.
The reduction of a K\"ahler manifold  is well-known (\cite{HLS}).
The reduction of an almost-K\"ahler  manifold is parallel and
even simpler since we do not need to check the integrability of $J$,
thus omitting the concrete construction of the Levi-Civita
connection.

Assume $(M,\omega,J, g)$ is an almost-K\"ahler  manifold with a free
$S^1$-Hamiltonian action also preserving $J$ (hence preserving $g$
as well). Let $X$ be the generating vector field of the
$S^1-$action. Suppose $\mu$ is a moment map. Then
$$\nabla \mu =JX.$$
It follows that each tangent space of $\mu^{-1}(0)/S^1$ has a
quotient complex vector space structure. Therefore $\mu^{-1}(0)/S^1$
has an induced almost- K\"ahler structure.
Let $(M, \omega, J)$ be an almost-K\"ahler  manifold with a
Hamiltonian $S^1$-action $\tau:S^1\times M\rightarrow M$, $\phi$ the
corresponding moment map, $\CC$ with its standard K\"ahler
structure.  Define an $S^1$-action

\begin{equation}  \begin{array}{llll}\tilde{\tau}:& S^1\times(M\times
\CC)&\rightarrow& M\times\CC \cr &(e^{i\theta},(x,z))&\rightarrow&
(e^{i\theta}x,e^{i\theta}z)\end{array}\end{equation}

\nono with moment map $\Psi(x,z)=\phi(x)+||z||^2$.
\defone\label{DefofKaehlercut}
  We call the
almost-K\"ahler manifold $\Psi^{-1}(\lambda)/S^1$ the
\textit{almost-K\"ahler  cut of $(M,\omega,J)$ along level set
$\lambda$}, denoted as $M_\lambda$. \deftwo

\nono Since the tangent space of $D_{\lambda}=\{(x,z)\in
M_\lambda|\phi(x)=\lambda\}/S^1$ is naturally identified with that
of reduced space $\phi^{-1}(\lambda)\slash S^1$, $D_{\lambda}$ is an
almost complex hypersurface of $M_{\lambda}$. Denote \oneeq
M'_\lambda=M_\lambda-D_{\lambda}. \twoeq Consider also the open
submanifold of $M$,
 \onedis M_0^\lambda=\{x\in M|\phi(x)<\lambda\}\twodis


$M'_\lambda$ and $M_0^\lambda$ come with induced almost complex
structure in the obvious way.  For example, the former
with the quotient almost-K\"ahler  structure, and the latter as a
subset of $M$.  The relation between symplectic structures of these
manifolds is clear.  As pointed out by E. Lerman \cite{Lerman},

\begin{equation} \label{canonicalsymplectomorphism}\begin{array}{llll}
T^{-1}:& M_0^\lambda&\rightarrow&
 M'_\lambda\cr
&x&\mapsto &(x,\sqrt{\lambda-\phi(x)})\end{array}\end{equation}

\nono is actually a symplectomorphism.
 However, $T$ is not an almost complex
isomorphism between the two manifolds with the almost complex structures
specified above.

Suppose the $S^1-$action on $(M, \omega, J)$ extends to an almost
complex action of $\mathbb C^*$ still denoted by $\tau$. For $z\in\CC^*$,
let $\tau_z$ be the corresponding almost complex isomorphism of $M$.  As observed in
\cite{BGL}, $\tau_{e^t}$ is generated by $\nabla\phi$.  Let

\oneeq\label{Prolongingequation0} \iota_0: M_0^\lambda\rightarrow M,\hskip 2mm
\iota_0(p)=\tau^{-1}_{\sqrt{\lambda-\phi(p)}}(p)\twoeq

\nono The central result in \cite{BGL}
we are going to use is the following:

\thmone\label{Kaehlercut}\textup{(\cite{BGL}, Theorem 2.1) With
notations above, then}

\oneeq \iota_0T:M'_\lambda\rightarrow \iota_0(M_0^\lambda)\twoeq

\nono\textup{is an almost complex isomorphism. }\thmtwo

Roughly speaking, this theorem asserts that
an open part of the almost-K\"ahler cut is biholomorphic to the
open set of $M$ which is the union of $M_0^\lambda$ and the unstable manifold
starting therein.  For readers' convenience, we include a sketch of the proof in \cite{BGL} with adaption to
the almost-K\"ahler case.
Consider the biholomorphic map
$$f:M\times \CC^*\rightarrow M\times \CC^*$$

\nono defined by $f(p,z_0)=(\tau_{z_0}p,z_0)$.  Consider actions $R$ and $D$ of $\CC^*$
on $M\times\CC^*$ defined as
$$R(z,(p,z_0))=(p,zz_0),$$
\nono and the diagonal action:
$$D(z,(p,z_0))=(\tau_zp,zz_0),$$

\nono then $f$ intertwines $R$ and $D$, i.e. $f\circ R(z,-)=D(z,f(-))$.
  Therefore, the pullback by $f$ of the natural product symplectic structure is invariant under the $S^1$-action on the $\CC^*$
factor (coming from the $\CC^*$-action restricted to the unit circle) with moment map
$$\tilde{\Psi}(p,z)=\phi(\tau_zp)+|z|^2$$

\nono Here $f$ maps the level set of $\tilde{\Psi}^{-1}(\lambda)$ onto the level set $\Psi^{-1}(\lambda)$ and induces an almost-
K\"ahler isomorphism:
$$h:\tilde{\Psi}^{-1}(\lambda)/S^1\rightarrow \Psi^{-1}(\lambda)/S^1$$

Note that the right-hand-side is identified with the almost-K\"ahler cut as desired.  Now we identify the left-hand-side
 of our isomorphism with the original manifold and
$h$ with $(\iota_0\circ T)^{-1}$ as follows:  remember we removed the origin of $\CC$ in the product manifold,
therefore, the $S^1$-quotients are actually free without cut-divisor.  Thus we have the following identification, where we identify the quotient of
$\CC^*$ by $S^1$ with the positive real ray:

$$\tilde{\Psi}^{-1}(\lambda)/S^1=\{(p,e^t)|\phi(\tau_{e^t}p)+e^{2t}=\lambda\}$$

\nono and,

$$\Psi^{-1}(\lambda)/S^1=\{(p,e^t)|\phi(p)+e^{2t}=\lambda\}$$

\nono For $\tilde{\Psi}^{-1}(\lambda)/S^1$, notice that as an almost-complex manifold, it is
$M\times\CC^*//\CC^*$, identified with $\iota_0(M_0^\lambda)$ with restricted almost-complex
structure from $M$ (since $f$ is a biholomorphism, what we did on the complex structure is simply multiplying $M$ by a
$\CC^*$ factor and then quotienting it out).  Under these identifications and (\ref{canonicalsymplectomorphism}),
$h(p,e^t)=(\tau_{e^t}p,e^t)$, hence

$$h^{-1}(p,e^t)=(\tau^{-1}_{\sqrt{\lambda-\phi(p)}}(p),\sqrt{\lambda-\phi(p)})$$

Hence $h=(\iota_0\circ T)^{-1}$, and in particular, the domain of $h$ is $\iota_0(M_0^\lambda)$.  Note that $T$
is just a step of identification which drops or picks up the second factor.
This proves the theorem.  For more details one is referred to \cite{BGL}.

\begin{remark}
We adopted different notations in Theorem \ref{Kaehlercut} from
what is in \cite{BGL}. In this note we use the notation
$M'_{\lambda}$ in place of the almost- K\"ahler  manifold
``$M_0^{\lambda}$ with its $M^{\lambda}$ complex structure" just for
convenience.  The content of the statement is completely the same as
it is in \cite{BGL}.
\end{remark}

\begin{example} The simplest example of K\"ahler cut is to consider the
symplectic cut of $\CC$ with standard K\"ahler structure at
$||z||=r$.  This induces a symplectomorphism from $\{||z||<r\}$ to
$\CC P^1\backslash \{pt\}$ with symplectic form
$\frac{r^2}{\pi}\omega_{FS}$.  On the other hand, sterographic
projection gives a biholomorphic map between $\CC$ and $\CC
P^1\backslash \{pt\}$, while these structures give a K\"ahler
structure in the cut space. Theorem \ref{Kaehlercut} gives a
formula between the moment map of $S^1$-action and the
biholomorphism, which could be viewed as a generalization of
stereoprojection.  Note that the same holds for $\CC^n$, and after
the cut at $||z||=r$, the symplectic area of the generator of
$H^2(\CC P^n)$ is $r^{2}$.
\end{example}

\subsection{A variation} \nono
In this subsection, $(M, \omega)$ is a symplectic manifold with a
preferred metric $g$.

\defone\textup (Extendable $S^1$-action) Suppose $J$ is an almost complex structure on $(M, \omega, g)$. A
  Hamiltonian $S^1$-action
$\tau$ on $M$  with a moment map $\phi$ is called $J-$extendable if
it is extended to an almost complex  semi-group action of the
complement of open unit disk in $\CC$, $\tilde{\tau}:\CC\backslash
D_1\times M\rightarrow M$, such that the gradient flow of the moment
map $\phi$ generates the infinitesimal action of
$\tilde{\tau}_{e^t}$. Here $S^1$ is identified with the unit circle
of $\CC$.
\deftwo


$(M,\omega, J)$ is called \textit{asymptotically almost-K\"ahler  (AAK)}
if there is a compact subset
$E$, where $(M\backslash E, \omega, J)$ is an almost--K\"ahler
manifold.
To avoid possible confusion, an extendable $S^1$-action on an AAK
manifold always acts only on $M\backslash E$, and the reference
metric is taken to be the one induced by $(\omega,J)$. It is easy to
see that asymptotically standard manifolds satisfies all assumptions
above: (see also \cite{BGL})


\lemmaone\label{symplectic form with S^1 symmetry} An asymptotically
standard manifold is AAK and the standard rotation is an extendable
$S^1$-action. \lemmatwo

\begin{remark}\label{more}\textit{
This lemma  is still valid if $(M,  \omega, J)$  is asymptotically
of the form  $(\RR^{2n}, \omega=-\sqrt{-1}\pl\bpl F(||z||^2) ,
J_0)$,
 where
 $F$ is a real valued function such that
$\frac{1}{\epsilon}>F'(x)+xF''(x)\geq\epsilon>0$ is bounded.
}\end{remark}

For our purpose, we need an adapted version of Theorem \ref{Kaehlercut} by allowing an extendable $S^1$-action on only the
non-compact part of the AAK manifold.  We first fix some notations for the
 convenience of exposition.  In the rest of the note, we make the convention that
$=_d,=_s,=_h,=_k$ denotes diffeomorphism, symplectomorphism,
biholomorphism and almost-K\"ahlerian isomorphism, respectively.
We also use the following notation:

\oneeq M_I=\{x\in M| \phi(x)\in I\subset\RR\}\twoeq

\nono With the above preparation, we are ready to state our lemma:

\lemmaone\label{KeyLemma} Let $(M,\omega, J)$ be an AAK manifold
with an extendable $S^1$-action on the complement of a compact set
$\tau: S^1\times M\backslash K\rightarrow M\backslash K$.  Then
there is an $\omega_{red}$-tamed almost complex structure $J_{red}$
on $M_\lambda$, such that
\[(M'_\lambda, J_{red})=_h(M, J)\]
\lemmatwo

\pfone  Let $\phi:M\backslash K\rightarrow \RR$ be the moment map on
the asymptotic part of $M$.  Suppose the $S^1$-action is extendable
on $\phi\geq\lambda_0$, choose $\lambda\gg\lambda_0$, let $\Psi:
M_{\geq\lambda_0}\times \CC\rightarrow \RR$ be the moment map as in
the definition of K\"ahler cut in section 2. Decompose

\onedis \Psi^{-1}(\lambda)=\{(x,z)|\phi(x)+|z|^2=\lambda\}\twodis
\oneeq=\{(x,z)|\phi(x)\leq\lambda-1\}\cup\\\{(x,z)|\lambda>\phi(x)>\lambda-1\}\cup\{(x,z)|\phi(x)=\lambda\}\twoeq
\onedis=A_1\cup A_2\cup A_3\twodis

\nono and $\Psi^{-1}(\lambda)/S^1$ has a natural symplectic form
$\omega_{red}$ from symplectic reduction.

Note that $\iota_0$ itself is identity on $\phi=\lambda-1$, to
define an $\omega_{red}$-tamed almost complex structure on
$M_\lambda$ , we first induce an almost complex structure
$T^{-1}_*(J)$ on $A_1\slash S^1$ by restricting the
symplectomorphism (\ref{canonicalsymplectomorphism}) to level set
below $\lambda-1$.  For the non-compact part, we denote
$N_{[\lambda-1,\lambda]}$ to be the reduction of

\oneeq M_{[\lambda-1,\lambda]}=\{x\in M\backslash
K|\lambda-1\leq\phi(x)\leq\lambda\}\twoeq

\nono Moreover, by the assumption of $\lambda$, the reduction is
almost-K\"ahlerian.  In particular, by Theorem \ref{Kaehlercut} and
(\ref{canonicalsymplectomorphism}), the reduced space would have
the following property:

\oneeq\label{holomorphism} (N_{[\lambda-1,\lambda)}, J_{red})=_h
(M_{[\lambda-1,\infty)}, J)\twoeq

\oneeq\label{symplectomorphism}
(N_{[\lambda-1,\lambda)},\omega_{red})=_s (M_{[\lambda-1,\lambda)},
\omega)\twoeq

\nono where the biholomorphism in (\ref{holomorphism}) is induced by
$(\iota_0 T)^{-1}|_{\RR^{2n}_{[\lambda-1,\infty)}}$ and the
symplectomorphism in (\ref{symplectomorphism}) by
$T^{-1}|_{M_{[\lambda-1,\lambda)}}$. Since the reduced space is an
almost-K\"ahler  manifold, \textit{a priori} we know the pair on the
left hand side of both equations above will consist of an
almost-K\"ahler pair.  Therefore, we have almost-K\"ahler
isomorphism via $T^{-1}$:

\[(M_{[\lambda-1,\lambda)},\omega,(\iota_0)^{-1}_*(J))=_k(N_{[\lambda-1,\lambda)},\omega_{red},J_{red})\]

In sum, we induce an $\omega_{red}$-tamed almost complex structure
$\tilde{J}$ on $(M_\lambda,\omega_{red})$ using the $C^1$-diffeomorphism:

\begin{equation}\label{defining complex structure} F=\iota_0\# id:M_{<\lambda}\rightarrow M\end{equation}
\[ \iota_0\# id= \left\{ \begin{aligned}
         \iota_0, &\text{         }\phi(x)\geq\lambda-1 \\
         id, &\text{         }\phi(x)\leq\lambda-1
                          \end{aligned} \right.
                          \]

\nono  To get a smooth almost-complex structure, one only
needs to perturb $F$ and use the openness of uniformly tamed
almost-complex structure with taming constant less than a fixed
$\beta>0$. In the following such a choice of perturbation would
not matter so we don't need to specify our choice of
perturbation.  Thus
$\tilde{J}$ is the tamed almost-complex structure as required.\pftwo

\rkone\label{tamingconstantofcut} This diffeomorphism is in fact the composition of $T$ and a time-$1$
map of an integral flow of a continuous vector field on $M$, see \cite{BGL} for a flow representation of
$\iota_0$.  Note that in the special case when $M$ is asymptotically
 standard, the taming and bounding constant are unchanged when almost-K\"ahler
 cut is performed along a standard sphere $M_\lambda\subset \RR^{2n}$.
  The reason is that the pair are unchanged in the compact set $M_{<\lambda}$;
 outside it is a subset of Euclidean space before the cut and complex projective space
 after the cut, with all corresponding structures, therefore the taming constant remains to be $1$.
\rktwo


\section{Review of Bangert's Results}

In this section we collect definitions and facts in \cite{Bangert}
for the proof of Theorem \ref{main}, and we refer the reader to
proofs therein for details.  We adopt standard notations in
\cite{Federer}.  For example, if $S$ is an $m-$current in $(M,g)$, $\MM(S)$
denotes its mass.  When $S$ is rectifiable, it is representable by integration and can be written as
$S=||S||\wedge\overrightarrow{S}$.  Here $||S||$ is a canonically associated Borel measure on $M$ finite on compact sets,
 and $\overrightarrow{S}:M\rightarrow \bigwedge^m M$ is the unit vector field orienting the
associated set of $S$.  If we have moreover a Borel subset $B$ on $M$, the \textit{interior multiplication}
 $S\llcorner B$ is the unique current represented as $S\llcorner B=\chi_B||S||\wedge\overrightarrow{S}$,
where $\chi_B$ is the characteristic function of $B$.  For a proper Lipschitz map $f$ and a rectifiable
current $S$, the push-forward current is denoted as $f_\#(S)$, which is also rectifiable.
See \cite{Federer} or \cite{Morgan} for the precise definitions.  Next we recall the notion of
quasi-minimality:

\begin{definition}\textup{ Let $(M, g)$ be a Riemannian manifold and $Q>1$. A rectifiable
current $S$ in $M$ is called quasiminimizing with constant $Q$, or
simply $Q$-minimizing if}
\[\MM(S\llcorner B)\leq Q\MM(X)\]
{when $B$ is a Borel subset of $M$ and $X$ a rectifiable current in
$M$ satisfying $\partial X=\partial(S\llcorner B)$}
\end{definition}

\begin{lemma}\label{J-curvesareQ-minimizing}
Let $\omega$ be an exact symplectic form, J an almost complex
structure on $(M,g)$. Suppose they form a uniformly tamed pair. Let
$(S,j)$ be a compact Riemann surface with orientation induced by
complex structure $j$, and $f:(S,j)\rightarrow (M,J)$ a
pseudo-holomorphic map.  Then the rectifiable current
$f_{\#}(S)$ is quasi-minimizing with
$Q=\alpha\beta/C$, with the notations in Definition \ref{Uniformlytamedpair}.
\end{lemma}

\nono Now a very useful lemma reads:

\lemmaone\label{Constraintlemma} For all real constants $Q>1, t>0$,
and integers $k\geq 1$, there exist $c=c(Q,t,k)\in (0,1)$ such that
the following is true for all $R>0$:  If $S$ is a $Q$-minimizing
rectifiable k-current in Euclidean $\RR^n$ s.t.

\onedis supp(\partial S)\subseteq B(cR), M(S)\leq tR^k,\twodis

\nono then

\onedis supp(S)\subseteq B(R).\twodis \qed\lemmatwo

We could adapt Lemma \ref{Constraintlemma} to asymptotically standard
manifold $M$ for $J$-holomorphic curves in our case.  For convenience, in the rest of the paper we denote
the closure of $(\phi_M^{-1}(\RR^{2n}\backslash
B_{R}))^c\subset M$ by $V_R$, $R>\bar{R}$ (See Definition
\ref{noncompact} for notations).

Suppose $J$ is a given uniformly tamed almost-complex structure, then $J$ restricts
to the non-compact standard part of $M$, therefore extends to give an
almost-complex structure $J'$ of $(\RR^{2n}, \omega_0)$ via
$\phi_M$.  By Lemma \ref{uniformlytamenessandcompactness} and
contractibility of $\mathscr{J}(\omega_0)$, we could choose such
extension by contracting $\mathscr{J}_{\beta, c}$ to $J_{std}$, so
that the almost-complex structure on the sphere $\partial V_R$
extends to the ball $B_R$.  Since such contraction can be chosen
once and for all and the trace is compact, we could assume the
extended $J'$ is uniformly tamed with taming and bounding constants
independent of $R$. Moreover, we thus obtain a constant
$c'(Q,J',t)>0$ from Lemma \ref{Constraintlemma}, such that for any
$J'$-holomorphic curves $S'$ with $supp(\partial S')\subseteq B(c'R),
M(S')\leq tR^2$, one has $S'\subset B_R\subset \RR^{2n}$.  Now for
any $R>\frac{1}{c'}\bar{R}$, let $S$ be a $J$-holomorphic curve with
boundary in $V_{c'R}$, then $S$ is
 divided into $S\cap V_{c'R}$ and $S\backslash V_{c'R}$.  We could apply Lemma \ref{J-curvesareQ-minimizing} and Lemma
 \ref{Constraintlemma} to $S\backslash V_{c'R}$ by identifying $S\backslash V_{c'R}$ with $J$-holomorphic
curve with boundary on $V_{c'R}$ in $\RR^{2n}$.  Therefore,

\lemmaone\label{GeneralizedConstraintlemma} Lemma \ref{Constraintlemma} is true for $J$-holomorphic curves in asymptotically standard
manifolds $(M, \omega, g)$ when $J$ is given, uniformly tamed by
$\omega$, and $c'R>\bar{R}$. \lemmatwo

\rkone\label{tamemanifold} Lemma \ref{GeneralizedConstraintlemma}
in particular shows an asymptotic standard manifold $M$ with a uniformly
tamed almost complex structure is a \textit{tame almost complex
manifold}, in the sense of Sikorav \cite{S}, which asserts that for
any compact set $K\subset M$ and positive class $C$, there is
another compact subset $K'\subset M$, such that every
$J$-holomorphic curve in the class $C$ intersecting $K$ is contained in
$K'$.  See \cite{Bio} for another generalization of Lemma
\ref{Constraintlemma} in Stein manifolds.\rktwo

\nono To deal with non-compact target space, Bangert proved the
following limiting process holds:

\begin{prop}\label{disklemma}
Let $(M, g, J)$ be an almost-complex manifold with a J-invariant
metric $g$. Suppose there is a sequence $R_j\rightarrow\infty$,
constant $C>0$, $J$-holomorphic maps $h_j:D_j\rightarrow M$ defined
on topological disks $D_j\subset \CC$, where $h_j(0)$ lies in a
compact subset $K\subset M$, such that
\begin{enumerate}
\item $dist_g(h_j(0),h_j(\partial D_j))\geq \alpha R_j$  and
\item $area_{g_0}(h_j)\leq CR_j^2$
\end{enumerate}

\nono Then $(M,g,J)$ is not uniformly hyperbolic.
\end{prop}

\nono With the aid of this proposition, the problem is reduced to
finding such $D_j$'s.  Note that Bangert did not state this
proposition in full generality, but the proof is easily seen to be
valid (cf. \cite{Bangert}, pp. 39, Proof of Proposition 2.7).

\section{Proof of Theorem \ref{main}}

\pfone  We follow Bangert's idea to construct holomorphic disks as
is required
 in Proposition \ref{disklemma}.  Let $J_{std}$ denote
the standard complex structure on $\RR^{2n}$ or $\CC P^n$, depending
on the space
 we are talking about.  We also use notations in Section 4.

Now for any $R>\bar{R}$, and fixed $J$ uniformly tamed by $\omega$
in the metric $g_0$, with bounding constant $\beta$ and taming
constant $c$, we can find $J_R$, such that the following holds:

\begin{enumerate}[(1)]
\item $J_R(x)=J_{std}(x)$, if $x\notin V_{2R}$ and $J_R=J$ if $x\in
V_R$.
\item $J_R$ is uniformly bounded and uniformly tamed by $\omega$
with taming constant $\beta'(\beta, c)$ and bounding constant
$c'(\beta, c)$ independent of $R$.
\end{enumerate}

This catenation is essentially proved in \cite{Bangert}.  Indeed, we
have a map $T:N(\partial V_R)\cup \RR^{2n}\backslash
B_{2R}\rightarrow \JJ(\omega_0)$, with
 $T(N(\partial V_R))\subset \Jst$ and $T(\RR^{2n}\backslash B_{2R})=J_{std}$, where $N(\partial V_R)$ is a neighborhood of $\partial V_R$.
As in the proof of Lemma \ref{GeneralizedConstraintlemma}, we
extend this map smoothly by contraction of $\Jst$ to $J_{std}$ to
obtain $J_R$.  Again since the catenation is constrained in a
compact subset, that is, the trace of the contraction, $J_R$ is a
uniformly tamed almost-complex structure with taming and bounding
constants independent of $R$.

For a fixed $J_R$ as above, we could perform almost-K\"ahler cut at
$\partial V_{3R}$, resulting in a new almost-K\"ahler manifold
$(M'_R,\omega_{red}, J_{red})$, which is almost-complex isomorphic
to $(M, J_R)$ by Lemma \ref{KeyLemma}, and preserves the taming
constant by Remark \ref{tamingconstantofcut}.
  It is readily seen that a
neighborhood of cut-divisor in $M_R$ is identified with
 a neighborhood of a divisor $\CC P^{n-1}\subset(\CC
P^n, \omega_{FS} ,J_{std})$ as an open K\"ahler manifold.  Therefore
the two-point Gromov-Witten invariant of $M_R$ in the class $t\in
H_2(M_R,\ZZ)$, the dual of the cut-divisor, is non-zero by Theorem
5.1.2 in \cite{HR}. Thus there is a stable curve $C$ passing
through a fixed point $x\in K$ in class $t$.  Consider the
components $\{C_i\}_{i=1}^l$ of $C$.

\begin{claim} There is exactly one component of $C$ which intersects the cut-divisor, and the
intersection is a unique point.  Moreover, such a component
intersects $K$ as well.
\end{claim}
\pfone

If a component $C_i$ intersects the cut-divisor, it either
intersects positively with or lies entirely in the cut-divisor $(\CC
P^{n-1}, J_{std})$. In the latter case, the intersection number of
$C_i$ and the cut divisor is also positive as the normal bundle of
the cut divisor is positive. Since the intersection number of $C$
and the cut divisor is $1$, the first two assertions follow.

The last assertion is proved if we can show that it is impossible
for any component to lie entirely in $M_R\backslash K$. This is
because that $M_R\backslash K$ is a disk bundle over the cut
divisor,  the class of any  component  in $M_R\backslash K$ must be
a positive multiple of $t$.  Since  $C$ also passes through $x\in
K$, there must exist other non-constant components. But these
components have non-positive total symplectic area, which is
impossible.
 \pftwo

We consider this unique component intersecting the cut-divisor.  This
holomorphic sphere gives a $J_R$-holomorphic line $f_R$ in $M$ with
$f_R(0)\in K$ and $f_R(\infty)$ the intersection with cut-divisor by
a reparametrization. Let $L_R:=f_R(\CC)\cap
\partial V_{R/2}$, and by Sard's theorem, we may assume
$f_{R}^{-1}(L_{R})$ is a collection of embedded circles up to an
arbitrary small change of $R$.  We choose the largest topological
disk $0\in D_R\sq \CC$ with boundary on $f_{R}^{-1}(L_{R})$. It is
clear from the claim above that the energy, or equivalently, the
symplectic area of $S_R:=(f_{R})_\#([D_R])$ is at most equal to that
of the cut-divisor, which is $\tilde{\nu} R^2$, $\tilde{\nu}$ an
absolute constant.  Note that the area of $S_R$ is bounded by
 the $\tilde{g}$-area induced by $(\omega, J)$, up to a constant $\bar{\nu}(\beta,
 C)$.  Therefore Remark \ref{Theremark} yields:

\oneeq\label{quasiminimal} {\bf M}(S_R)\leq \nu R^2 \twoeq \nono
where $\nu=\nu(\beta, C)$.

The general procedure above does not give holomorphic disks with
respect to the almost complex structure $J$ because $S_R$ might
touch the ``shell" $V_{2R}\backslash V_{R}$ where the almost complex
structure $J_{R}$ does not coincide with $J$.  To get a genuine
$J$-holomorphic disk, consider a subset of $S_R$.  We have the
uniform quasi-minimality constant $Q$ and $\nu$ satisfying
(\ref{quasiminimal}), Lemma \ref{GeneralizedConstraintlemma} thus
gives a uniform constant $c\in(0,1)$, such that when $f(\partial
D')\in V_{cR}$, we have $f(D')\in V_R$.  Hence by considering the
intersection $\partial V(cR)\cap S_R$, and repeating the above
argument with possibly one more application of Sard's theorem, we
have a holomorphic disk $\bar{D}_R$ lying entirely in $V_R$,
therefore a $J$-holomorphic disk, with boundary on $V_{cR}$.  It
follows then from Proposition \ref{disklemma} that $(M,g,J)$ is not
uniformly hyperbolic.

 \pftwo

\rkone We in fact did not use Lemma \ref{KeyLemma} in an essential
way in the proof above.  Instead, one could just use stereographic
projection to assign complex structure on the cut space to argue for
uniform tameness. However, our proof is generalized in a
straightforward way to asymptotic $\RR^{2n}$ with standard metric
and the symplectic form in Remark \ref{more}.  One easily see that
the neighborhood of cut-divisor as a complex manifold is still
$\mathscr{O}(1)$ therefore positive.  Moreover, we have also the
estimate (\ref{quasiminimal}).  This is equivalent to an estimate
of the symplectic area, thus the associated $\tilde{g}$-area of the
cut-divisor. We could view the cut-divisor as $\CC P^{n-1}$ with
corresponding symplectic form $\omega$, and by a homology argument
we reduce the problem to the estimate of the line class pairing with
$\omega$. Therefore we could just take the holomorphic disk in
$x_1x_2$-plane, and by assumption on $F$, the associated
$\tilde{g}=(F'+xF'')\cdot g_{eud}$ is bounded and equivalent to
$g_{eud}$, as desired.\rktwo

\end{document}